\documentclass{article}
\usepackage{amsfonts,amsmath,graphicx,amsthm,verbatim,mathrsfs}
\usepackage[dvips]{color}
\newtheorem{defn}{Definition}[section]
\newtheorem{conj}[defn]{Conjecture}
\newtheorem{prop}[defn]{Proposition}
\newtheorem{lemma}[defn]{Lemma}
\newtheorem{thm}[defn]{Theorem}
\newtheorem{cor}[defn]{Corollary}

\numberwithin{equation}{section}
\newcommand{\com}[2]{\left[\,#1\,,#2\,\right]}
\newcommand{\inprod}[2]{\left(\,#1\,,#2\,\right)}
\DeclareMathOperator{\sgn}{sgn}
 
\DeclareMathOperator{\prob}{P} \DeclareMathOperator{\airy}{Ai}
\newcommand{\rchi}{\raisebox{.4ex}{$\chi$}}

\begin{document}
\title{Distribution Functions for Edge Eigenvalues in Orthogonal and
Symplectic Ensembles: Painlev\'e Representations}
\date{}
\author{\textsc{Momar Dieng} \\ \\ \it{Department of Mathematics} \\
\it{University of California, Davis, CA 95616, USA}
\\ \it{E-mail address: momar@math.ucdavis.edu}}
\maketitle

\begin{abstract}
We derive Painlev\'e--type expressions for the distribution of the
$m^{th}$ largest eigenvalue in the Gaussian Orthogonal and
Symplectic Ensembles in the edge scaling limit. The work of
Johnstone and Soshnikov (see \cite{John1}, \cite{Sosh2}) implies
the immediate relevance of our formulas for the $m^{th}$ largest
eigenvalue of the appropriate Wishart distribution.

\end{abstract}

\section{Introduction}
\begin{figure}[htbp]
\includegraphics{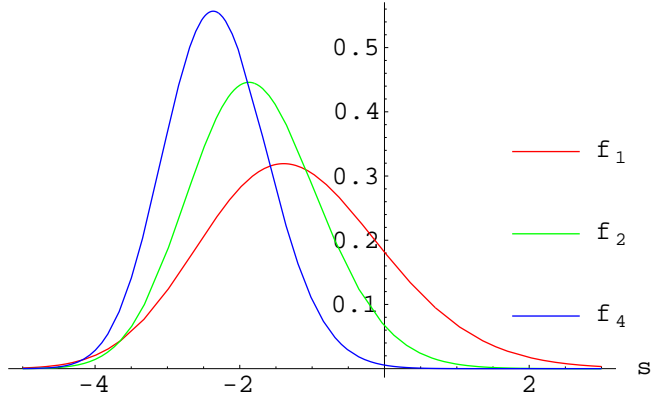}
\caption{Tracy--Widom Density Functions} \label{TWdensities}
\end{figure}

The Gaussian $\beta$--ensembles are probability spaces on
$n$-tuples of random variables $\{\lambda_{1},\ldots,
\lambda_{N}\}$, with joint density functions
\begin{equation} P_{N\,\beta}(\lambda_{1},\ldots,\lambda_{N})= P_{N\,\beta}(\vec{\lambda}) = C_{N\beta}\,\exp\left[-\frac{1}{2}\beta\,\sum_{j=1}^{N}\lambda_{j}^{2}\right]\prod_{j<k}|\lambda_{j}-\lambda_{k}|^{\beta}.\label{jointdensity}\end{equation}
The $C_{N\beta}$ are normalization constants, and by setting
$\beta=1, 2, 4$ we recover the  \emph{Gaussian Orthogonal
Ensemble} ($\textrm{GOE}_{N}$), \emph{Gaussian Unitary Ensemble}
($\textrm{GUE}_{N}$), and \emph{Gaussian Symplectic Ensemble}
($\textrm{GSE}_{N}$), respectively. We restrict ourselves to those
three cases in this paper, and refer the reader to \cite{Dumi1}
for recent results on the general $\beta$ case. Originally the
$\lambda_{j}$ are eigenvalues of randomly chosen matrices from
corresponding matrix ensembles, so we will henceforth refer to
them as eigenvalues. With the eigenvalues ordered so that
$\lambda_{j}\geq\lambda_{j+1}$, define
\begin{equation}\hat{\lambda}_{m}^{(N)}=\frac{\lambda_{m}-\sqrt{2\,N}}{2^{-1/2}\,N^{-1/6}},\end{equation}
to be the rescaled $m^{th}$ eigenvalue measured from edge of
spectrum. A standard result of Random Matrix Theory about the
distribution of the largest eigenvalue in the $\beta$--ensembles
is that
\begin{equation} \hat{\lambda}_{1}^{(N)}\xrightarrow{\mathscr{D}}\hat{\lambda}_{1},\end{equation}
whose law is given by the Tracy--Widom distributions.
\begin{thm}[Tracy, Widom \cite{Trac3},\cite{Trac2}]
\begin{equation}F_{2}(s):=\prob_{_{\textrm{GUE}}}(\hat{\lambda}_{1}\leq s)=\exp\left[-\int_{s}^{\infty}(x-s)\,q^{2}(x)d\,x\right]
\label{guemax},\end{equation}
\begin{equation}F_{1}^{2}(s):=\left[\prob_{_{\textrm{GOE}}}(\hat{\lambda}_{1}\leq s)\right]^{2}=F_{2}\cdot\exp\left[-\int_{s}^{\infty}q(x)d\,x\right]\label{goemax},\end{equation}
\begin{equation}F_{4}^{2}(\frac{s}{\sqrt{2}}):=\left[\prob_{_{\textrm{GSE}}}(\hat{\lambda}_{1}\leq s)\right]^{2}=F_{2}\cdot\cosh^{2}\left[-\frac{1}{2}\int_{s}^{\infty}q(x)d\,x\right]
\label{gsemax}.\end{equation}
\end{thm}
The function $q$ is the unique (see \cite{Hasti1},\cite{clar1})
solution to the Painlev\'e II equation
\begin{equation} q'' = x\,q + 2\,q^{3},\label{pII}\end{equation}
such that $q(x)\sim \airy(x)$ as $x\to\infty$, where $\airy(x)$ is
the solution to the Airy equation which decays like
$\frac{1}{2}\,\pi^{-1/2}\,x^{-1/4}\,\exp\left(-\frac{2}{3}\,x^{3/2}\right)$
at $+\infty$. The density functions $f_{\beta}$ corresponding to
the $F_\beta$ are graphed in Figure
\ref{TWdensities}.\footnote{The square root of 2 in the argument
of $F_{4}$ reflects a normalization chosen in \eqref{jointdensity}
to agree with Mehta's original one. It can be removed by choosing
a different normalization.}

Let $F_{2}(s,m)$ denote the distribution for the $m^{th}$ largest
eigenvalue in GUE. Tracy and Widom showed in \cite{Trac3} that if
we define $F_{2}(s,0)\equiv 0$, then

\begin{equation}F_{2}(s,m+1) -  F_{2}(s,m) = \frac{(-1)^{m}}{m\,!}\frac{d^{m}}{d\,\lambda^{m}}\,D_{2}(s,\lambda)\big{\vert}_{\lambda=1}\,,\quad
m\geq 0, \label{lderiv} \end{equation} where
\begin{equation}D_{2}(s,\lambda)=\exp\left[-\int_{s}^{\infty}(x-s)\,q^{2}(x,\lambda)d\,x \right],\label{D2}\end{equation}
and $q(x,\lambda)$ is the solution to \eqref{pII} such that
$q(x,\lambda)\sim~\sqrt{\lambda}\,\airy(x)$ as $x\to\infty$. An
intermediate step leading to \eqref{D2} is to first show that
$D_{2}(s,\lambda)$ can be expressed as a Fredholm determinant
\begin{equation}D_{2}(s,\lambda)=\det(I-\lambda\,\mathcal K_{\airy}),\label{fred}\end{equation}
where $\mathcal K_{\airy}$ is the integral operator with kernel
\begin{equation}\mathcal K_{\airy}(x,y)=\frac{\airy(x)\airy'(y)-\airy'(x)\airy(y)}{x-y}\,.\end{equation}
In the $\beta=1,4$ cases a result similar to \eqref{fred} holds
with the difference that the operators in $D_{\beta}(s,\lambda)$
have matrix--valued kernels (see e.g. \cite{Trac5}). In fact, the
same combinatorial argument used to obtain the recurrence
\eqref{lderiv} in the $\beta=2$ case also works for the
$\beta=1,4$ cases, leading to
\begin{equation}
F_{\beta}(s,m+1) -  F_{\beta}(s,m) =
\frac{(-1)^{m}}{m\,!}\frac{d^{m}}{d\,\lambda^{m}}\,D_{\beta}^{1/2}(s,\lambda)\big{\vert}_{\lambda=1}\,,\quad
m\geq 0,\, \beta=1,4, \label{betaderiv}\end{equation} where
$F_{\beta}(s,0)\equiv 0$. Given the similarity in the arguments up
to this point and comparing \eqref{D2} to \eqref{guemax}, it is
natural to conjecture that $D_{\beta}(s,\lambda), \beta=1,4,$ can
be obtained simply by replacing $q(x)$ by $q(x,\lambda)$ in
\eqref{goemax} and \eqref{gsemax}. However the following
conjecture, which had long been in the literature, and whose
verification is the content of Corollary~\eqref{interlacingcor},
hints that this cannot be the case:
\begin{conj}[Baik, Rains \cite{Baik1}]
In the appropriate  scaling limit, the distribution of the largest
eigenvalue in GSE corresponds to that of the second largest in
GOE. More generally, the joint distribution of every second
eigenvalue in the GOE coincides with the joint distribution of all
the eigenvalues in the GSE, with an appropriate number of
eigenvalues. \label{baikconj}\end{conj} Forrester and Rains
subsequently proved (see \cite{Forr1}) the equivalence of
alternate GOE eigenvalues and GSE eigenvalues at \emph{finite} $N$
ensemble level, lending weight to Conjecture~\eqref{baikconj}.
This so--called ``interlacing property'' between GOE and GSE had
been noticed by Mehta and Dyson (see \cite{Meht3}).
Conjecture~\eqref{baikconj} does not agree with the formulae we
postulated for $D_{\beta}(s,\lambda), \beta=1,4$. Indeed,
combining the two leads to incorrect relationships between
derivatives of $q(x,\lambda)$ evaluated at $\lambda=1$. To be
precise, the conjecture is true for $D_{4}(s,\lambda)$ but it is
false for $D_{1}(s,\lambda)$. The correct forms for both
$D_{\beta}(s,\lambda), \beta=1,4$ are given below in
Theorem~\eqref{mainthm}.
\par This work also extends that of Johnstone in \cite{John1} (see also \cite{Elka1}), since $F_{1}(s,m)$
gives the asymptotic behavior of the $m^{th}$ largest eigenvalue
of a $p$ variate
 Wishart distribution on $n$ degrees of freedom with identity covariance. This holds under very
 general conditions on the underlying distribution of matrix entries by Soshnikov's
 universality theorem (see \cite{Sosh2} for a precise statement).
 In Table~\ref{table}, we compare our distributions to
 finite $n$ and $p$ empirical Wishart distributions as in \cite{John1}.

 \section{Statement of the Main Results}

\begin{thm}
The distributions for the $m^{th}$ largest eigenvalues in the
\textrm{GOE} and \textrm{GSE} satisfy the recurrence
\eqref{betaderiv} with
\begin{equation}D_{1}(s,\lambda)=D_{2}(s,\tilde{\lambda})\,\frac{\lambda - 1 -
\cosh{\mu(s,\tilde{\lambda})} +
\sqrt{\tilde{\lambda}}\,\sinh{\mu(s,\tilde{\lambda})}}{\lambda -
2}, \label{goedet}\end{equation}
\begin{equation}D_{4}(s,\lambda)=D_{2}(s,\lambda)\,\cosh^{2}\left(\frac{\mu(s,\lambda)}{2}\right),\label{gsedet}\end{equation}
where
\begin{equation}\mu(s,\lambda):=\int_{s}^{\infty}q(x,\lambda)d\,x, \qquad \tilde{\lambda}:=2\,\lambda-\lambda^{2},\end{equation}
and $q(x,\lambda)$ is the solution to \eqref{pII} such that
$q(x,\lambda)\sim~\sqrt{\lambda}\,\airy(x)$ as
$x\to\infty$.\label{mainthm}\end{thm}
\begin{cor}[Interlacing property]
\begin{equation}F_{4}(s,m)= F_{1}(s,2m),\quad m\geq 1.\end{equation}
\label{interlacingcor}\end{cor} In the next section we outline the
proof of these theorems. In the last, we present an efficient
numerical scheme to compute $F_{\beta}(s,m)$. We implemented this
scheme using MATLAB, and compared the results to simulated Wishart
distributions.
\begin{center}
\begin{figure}
  \includegraphics[height=65mm]{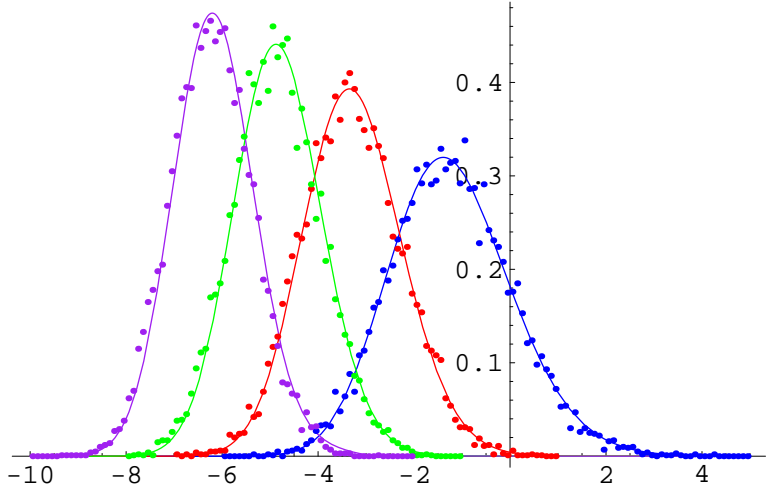}
\caption{$10^{4}$ realizations of $10^{3}\times 10^{3}$ GOE
matrices; the solid curves are, from right to left, the
theoretical limiting densities for the first through fourth
largest eigenvalue.}
\end{figure}
\end{center}

\section{Sketch of the Proofs}

\subsection{Distribution for the Next Largest Eigenvalues and\\
Finite $N$ Gaussian Ensembles} With the joint density function
defined as in \eqref{jointdensity}, let $J$ be an interval on the
real line, and $\rchi=\rchi_{_{J}}(x)$ its characteristic
function. We denote by $\tilde{\rchi}=1-\rchi$ the characteristic
function of the complement of $J$, and define
$\tilde{\rchi}_{_{\lambda}}=1-\lambda\,\rchi$. Furthermore, let
$E_{\beta,N}(m,J)$ equal the probability that exactly the $m$
largest eigenvalues of a matrix chosen at random from a (finite
$N$) $\beta$--ensemble lie in $J$. We also define

\begin{equation} E_{\beta,N}^{(\lambda)}(J)= \underset{x_{i}\in \mathbb
R}{\int\cdots\int} \tilde{\rchi}_{_{\lambda}}(x_{1})\cdots
\tilde{\rchi}_{_{\lambda}}(x_{N})\,P_{N\,\beta}(x_{1},\ldots,x_{N})\,d\,x_{1}\cdots
d\,x_{N}.\end{equation} For $\lambda=1$ this is just
$E_{\beta,n}(0,J)$, the probability that no eigenvalues lie in
$J$. The following propositions are easy combinatorial facts that
can be proved by induction (see e.g. \cite{Trac6}).

\begin{prop}
\begin{equation}
E_{\beta,N}^{(\lambda)}(J)=\sum_{k=0}^{N}(-\lambda)^{k}\binom{N}{k}\underset{x_{i}\in
J}{\int\cdots\int}
P_{N\,\beta}(x_{1},\ldots,x_{N})\,d\,x_{1}\cdots
d\,x_{N}.\end{equation}
\end{prop}

\begin{prop}
\begin{equation}\left.
E_{\beta,N}(m,J)=\frac{(-1)^{m}}{m\,!}\,\frac{d^{m}}{d\,\lambda^{m}}\,E_{\beta,N}^{(\lambda)}(J)\right|_{\lambda=1}\,
,\quad m\geq 0.\end{equation}
\end{prop}
The next step is to find a useful expression for the multiple
integral $E_{\beta,N}^{(\lambda)}(J)$. It turns out that through
standard RMT techniques (see e.g. \cite{Trac1}), the integral can
be expressed as the determinant of an operator on the Hilbert
space
 $L^{2}(J)\times~L^{2}(J)$. Let

\begin{equation}D_{\beta,N}(s,\lambda)=\det(I-\lambda\,\mathcal K_{\beta,N})\quad,\quad\beta=1,4,\label{fredholm}\end{equation}
for
\begin{equation} \mathcal K_{1,N} = \rchi\,\left(\begin{array}{cc}
S + \psi\otimes\epsilon\,\varphi & S\,D - \psi\otimes\varphi \\
\epsilon\,S - \epsilon + \epsilon\,\psi\otimes\epsilon\,\varphi &
S + \epsilon\,\varphi\otimes\psi
\end{array}\right)\,\rchi\label{goekernel}.\end{equation}
Here $\epsilon$ is the integral operator with kernel
$\epsilon(x-y)=\frac{1}{2}\sgn(x-y)$, and $D$ denotes the
differentiation operator $\frac{d}{d\,x}$, $S$ is the integral
operator with kernel
\begin{equation}S(x,y)=\frac{\varphi(x)\,\psi(y)-\psi(x)\,\varphi(y)}{x-y},\end{equation} and the functions
$\varphi$ and $\psi$ are

\begin{equation}\varphi(x)=\left(\frac{N}{2}\right)^{1/4}\,\varphi_{N}(x),\end{equation}

\begin{equation}\psi(x)=\left(\frac{N}{2}\right)^{1/4}\,\varphi_{N-1}(x),\end{equation}
where
\begin{equation}\varphi_{N}(x)=\frac{1}{\sqrt{2^{N}\,N!\,\sqrt{\pi}}}\,e^{-x^{2}/2}\,H_{N}(x),\end{equation}
and the $H_{N}(x)$ are the classical Hermite polynomials. This
implies that the $\varphi_{n}(x)$ are orthonormal with respect to
the Lebesgue measure on $\mathbb R$. Similarly, let
\begin{equation}\mathcal K_{4,N}=\frac{1}{2}\,\rchi\,\left(\begin{array}{cc} S + \psi\otimes\,\epsilon\,\varphi  & SD - \psi\otimes\,\varphi \\
 \epsilon\,S + \epsilon\,\psi\otimes\,\epsilon\,\varphi  &  S + \epsilon\,\varphi\otimes\,\epsilon\,\psi
 \end{array}\right)\,\rchi\label{gsekernel}.\end{equation}
 Following the same approach as in \cite{Trac1} and
\cite{Trac2}, we arrive at
\begin{equation}
E_{\beta,N}^{(\lambda)}(J)=D_{\beta,N}^{1/2}(s,\lambda).
\end{equation}

\subsection{Edge-Scaling}

\subsubsection{The GOE case: reduction of the determinant}
The above determinants are Fredholm determinants of operators on
$L^{2}(J)~\times~L^{2}(J)$. Our first task will be to rewrite
these determinants as those of operators on $L^{2}(J)$. This part
follows exactly the proof in \cite{Trac2}. To begin, note that
\begin{equation}\com{S}{D}=\varphi\otimes\psi +
\psi\otimes\varphi\label{sdcom}\end{equation} so that (using the
fact that $D\,\epsilon=\epsilon\,D=I$ )
\begin{eqnarray}
\com{\epsilon}{S} &=&\epsilon\,S-S\,\epsilon\nonumber\\
&=& \epsilon\,S\,D\,\epsilon-\epsilon\,D\,S\,\epsilon = \epsilon\,\com{S}{D}\,\epsilon\nonumber\\
 &=& \epsilon\,\varphi\otimes\psi\,\epsilon + \epsilon\,\psi\otimes\varphi\,\epsilon\nonumber\\
&=&\epsilon\,\varphi\otimes\epsilon^{t}\psi + \epsilon\,\psi\otimes\epsilon^{t}\,\varphi\nonumber\\
&=& - \epsilon\,\varphi\otimes\epsilon\,\psi -
\epsilon\,\psi\otimes\epsilon\,\varphi,
\label{escom}\end{eqnarray} where the last equality follows from
the fact that $\epsilon^{t}=-\epsilon$. We thus have

\begin{eqnarray*}
D\,\left(\epsilon\,S + \epsilon\,\psi\otimes\epsilon\varphi\right) & = & S + \psi\otimes\epsilon\varphi,\\
D\,\left(\epsilon\,S\,D - \epsilon\,\psi\otimes\varphi\right) & =
& S\,D - \psi\otimes\varphi.\end{eqnarray*} The expressions on the
right side are the top entries of $\mathcal K_{1,N}$. Thus the
first row of $\mathcal K_{1,N}$ is, as a vector,

\begin{equation*}D\,\left(\epsilon\,S + \epsilon\,\psi\otimes\epsilon\varphi,
\epsilon\,S\,D - \epsilon\,\psi\otimes\varphi\right).
\end{equation*} Now \eqref{escom} implies that

\begin{equation*}  \epsilon\,S + \epsilon\,\psi\otimes\epsilon\,\varphi =
S\,\epsilon -\epsilon\,\varphi\otimes\epsilon\,\psi.
\end{equation*} Similarly \eqref{sdcom} gives

\begin{equation*} \epsilon\,\com{S}{D} =
\epsilon\,\varphi\otimes\psi + \epsilon\psi\otimes\varphi,
\end{equation*} so that

\begin{equation*}\epsilon\,S\,D - \epsilon\,\psi\otimes\varphi = \epsilon\,D\,S
+ \epsilon\,\varphi\otimes\psi = S + \epsilon\,\varphi\otimes\psi.
\end{equation*}
Using these expressions we can rewrite the first row of $\mathcal
K_{1,N}$ as

\begin{equation*}D\,\left(S\,\epsilon - \epsilon\,\varphi\otimes\epsilon\psi, S
+ \epsilon\,\varphi\otimes\psi\right). \end{equation*} Applying
$\epsilon$ to this expression shows the second row of $\mathcal
K_{1,N}$ is given by

\begin{equation*}\left(\epsilon\,S - \epsilon +
\epsilon\,\psi\otimes\epsilon\varphi, S +
\epsilon\,\varphi\otimes\psi\right) \end{equation*} Now use
\eqref{escom} to show the second row of $\mathcal K_{1,N}$ is

\begin{equation*}\left(S\,\epsilon - \epsilon +
\epsilon\,\varphi\otimes\epsilon\psi, S +
\epsilon\,\varphi\otimes\psi\right). \end{equation*} Therefore,

\begin{eqnarray*}\mathcal K_{1,N} &=& \rchi\,\left(\begin{array}{cc}
D\,\left(S\,\epsilon - \epsilon\,\varphi\otimes\epsilon\psi\right) & D\,\left(S + \epsilon\,\varphi\otimes\psi\right) \\
 S\,\epsilon - \epsilon + \epsilon\,\varphi\otimes\epsilon\psi &  S + \epsilon\,\varphi\otimes\psi
\end{array}\right)\,\rchi \\
  & =& \left(\begin{array}{cc}
\rchi\,D & 0 \\
 0 & \rchi
\end{array}\right) \left(\begin{array}{cc}
\left(S\,\epsilon - \epsilon\,\varphi\otimes\epsilon\psi\right)\,\rchi & \left(S + \epsilon\,\varphi\otimes\psi\right)\,\rchi \\
 \left(S\,\epsilon - \epsilon + \epsilon\,\varphi\otimes\epsilon\psi\right)\rchi & \left(S + \epsilon\,\varphi\otimes\psi\right)\rchi
\end{array}\right).
\end{eqnarray*}
Since $\mathcal K_{1,N}$ is of the form $A\,B$, we can use the
fact that $\det(I-A\,B)=~\det(I-~B\,A)$ and deduce that
$D_{1,N}(s,\lambda)$ is unchanged if instead we take $\mathcal
K_{1,N}$ to be

\begin{eqnarray*} \mathcal K_{1,N} & = & \left(\begin{array}{cc}
\left(S\,\epsilon - \epsilon\,\varphi\otimes\epsilon\psi\right)\,\rchi & \left(S + \epsilon\,\varphi\otimes\psi\right)\,\rchi \\
 \left(S\,\epsilon - \epsilon + \epsilon\,\varphi\otimes\epsilon\psi\right)\,\rchi & \left(S + \epsilon\,\varphi\otimes\psi\right)\,\rchi
\end{array}\right) \left(\begin{array}{cc}
\rchi\,D & 0 \\
 0 & \rchi
\end{array}\right)\\
& = & \left(\begin{array}{cc}
\left(S\,\epsilon - \epsilon\,\varphi\otimes\epsilon\psi\right)\,\rchi\,D & \left(S + \epsilon\,\varphi\otimes\psi\right)\,\rchi \\
 \left(S\,\epsilon - \epsilon + \epsilon\,\varphi\otimes\epsilon\psi\right)\rchi\,D & \left(S + \epsilon\,\varphi\otimes\psi\right)\rchi
\end{array}\right).
\end{eqnarray*}
Therefore

\begin{equation}D_{1,N}(s,\lambda)=\det\left(\begin{array}{cc}
I - \left(S\,\epsilon - \epsilon\,\varphi\otimes\epsilon\psi\right)\,\lambda\,\rchi\,D & - \left(S + \epsilon\,\varphi\otimes\psi\right)\,\lambda\,\rchi \\
 - \left(S\,\epsilon - \epsilon + \epsilon\,\varphi\otimes\epsilon\psi\right)\,\lambda\,\rchi\,D & I - \left(S + \epsilon\,\varphi\otimes\psi\right)\,\lambda\,\rchi
\end{array}\right).\end{equation}
Now we perform row and column operations on the matrix to simplify
it, which do not change the Fredholm determinant. Justification of
these operations is given in \cite{Trac2}. We start by subtracting
row 1 from row 2 to get

\begin{equation*}  \left(\begin{array}{cc}
I - \left(S\,\epsilon - \epsilon\,\varphi\otimes\epsilon\psi\right)\,\lambda\,\rchi\,D & - \left(S + \epsilon\,\varphi\otimes\psi\right)\,\lambda\,\rchi \\
  - I + \epsilon\,\lambda\,\rchi\,D  & I \end{array}\right). \end{equation*}
 Next, adding column 2 to column 1 yields

 \begin{equation*} \left(\begin{array}{cc}
I - \left(S\,\epsilon - \epsilon\,\varphi\otimes\epsilon\psi\right)\,\lambda\,\rchi\,D - \left(S + \epsilon\,\varphi\otimes\psi\right)\,\lambda\,\rchi & - \left(S + \epsilon\,\varphi\otimes\psi\right)\,\lambda\,\rchi \\
  \epsilon\,\lambda\,\rchi\,D  & I \end{array} \right). \end{equation*}
Then right-multiply column 2 by $-\epsilon\,\lambda\,\rchi\,D$ and
add it to column 1 to get

\begin{equation*} \left(\begin{array}{cc}
I - \left(S\,\epsilon -
\epsilon\,\varphi\otimes\epsilon\psi\right)\,\lambda\,\rchi\,D +
\left(S + \epsilon\,\varphi\otimes\psi\right)\,\lambda\,\rchi\,\left(\epsilon\,\lambda\,\rchi\,D - I\right) & - \left(S + \epsilon\,\varphi\otimes\psi\right)\,\rchi \\
  0  & I \end{array} \right).\end{equation*}
Finally we multiply row 2 by $S + \epsilon\,\varphi\otimes\psi$
and add it to row 1 to arrive at

\begin{equation*} \det\left(\begin{array}{cc} I - \left(S\,\epsilon -
\epsilon\,\varphi\otimes\epsilon\psi\right)\,\lambda\,\rchi\,D +
 \left(S + \epsilon\,\varphi\otimes\psi\right)\,\lambda\,\rchi\,\left(\epsilon\,\lambda\,\rchi\,D - I\right) & 0 \\
  0  & I \end{array} \right). \end{equation*}
Thus the determinant we want equals the determinant of

  \begin{equation}I - \left(S\,\epsilon - \epsilon\,\varphi\otimes\epsilon\psi\right)\,\lambda\,\rchi\,D + \left(S +
  \epsilon\,\varphi\otimes\psi\right)\,\lambda\,\rchi\,\left(\epsilon\,\lambda\,\rchi\,D - I\right)
  \label{operator}.\end{equation}
So we have reduced the problem from the computation of the
Fredholm determinant of an operator on $L^{2}(J)~\times~L^{2}(J)$,
to that of an operator on $L^{2}(J)$.

\subsubsection{The GOE Case: differential equations}

Next we want to write the operator in \eqref{operator} in the form
\begin{equation} \left(I- \mathcal K_{2, N}\right)\left(I -
\sum_{i=1}^{L}\alpha_{i}\otimes\beta_{i}\right), \end{equation}
where the $\alpha_{i}$  and $\beta_{i}$ are functions in
$L^{2}(J)$. In other words, we want to rewrite the determinant for
the GOE case as a finite dimensional perturbation of the
corresponding GUE determinant. The Fredholm determinant of the
product is then the product of the determinants. The limiting form
for the GUE part is already known, and we can just focus on
finding a limiting form for the determinant of the finite
dimensional piece. It is here that the proof must be modified from
that in \cite{Trac2}. A little simplification of \eqref{operator}
yields
\begin{equation*}
I -
\lambda\,S\,\rchi-\lambda\,S\,\left(1-\lambda\,\rchi\right)\,\epsilon\,\rchi\,D
- \lambda\,\left(\epsilon\,\varphi\,\otimes\,\rchi\,\psi\right) -
\lambda\,\left(\epsilon\,\varphi\,\otimes\,\psi\right)\left(1-\lambda\,\rchi\right)\,\epsilon\,\rchi\,D.
\end{equation*}
Writing $\epsilon\,\com{\rchi}{D}+\rchi$ for $\epsilon\,\rchi\,D$
and simplifying $\left(1-\lambda\,\rchi\right)\,\rchi$ to
$\left(1-\lambda\right)\,\rchi$ gives

 \begin{align*}& I -
\lambda\,S\,\rchi - \lambda\,\left(1-\lambda\right)\,S\,\rchi -
\lambda\,\left(\epsilon\,\varphi\,\otimes\,\rchi\,\psi\right)
-\lambda\,\left(1-\lambda\right)\,\left(\epsilon\,\varphi\,\otimes\,\rchi\,\psi\right)
\\ & \qquad
-\lambda\,S\,\left(1-\lambda\,\rchi\right)\,\epsilon\,\com{\rchi}{D}
-
\lambda\,\left(\epsilon\,\varphi\,\otimes\,\psi\right)\,\left(1-\lambda\,\rchi\right)\,\epsilon\,\com{\rchi}{D}
 \\ & = I -  (2\lambda-\lambda^{2})\,S\,\rchi - (2\lambda-\lambda^{2})\,(\epsilon\,\varphi\,\otimes\,\rchi\,\psi)
  - \lambda\,S\,(1-\lambda\,\rchi)\,\epsilon\,\com{\rchi}{D} \\ &
  \qquad
  -  \lambda\,(\epsilon\,\varphi\,\otimes\,\psi)\,(1-\lambda\,\rchi)\,\epsilon\,\com{\rchi}{D}.\end{align*}

Define $\tilde{\lambda}=2\,\lambda-\lambda^{2}$ and let
$\sqrt{\tilde{\lambda}}\,\varphi\to\varphi$, and
$\sqrt{\tilde{\lambda}}\,\psi\to\psi$ so that
$\tilde{\lambda}\,S\to S$ and \eqref{operator} goes to
\begin{align*} I -  & S\,\rchi - (\epsilon\,\varphi\,\otimes\,\rchi\,\psi)
  - \frac{\lambda}{\tilde{\lambda}}\,S\,(1-\lambda\,\rchi)\,\epsilon\,\com{\rchi}{D} \\ & -
 \frac{\lambda}{\tilde{\lambda}}\,(\epsilon\,\varphi\,\otimes\,\psi)\,(1-\lambda\,\rchi)\,\epsilon\,\com{\rchi}{D}.\end{align*}
Now we define $R:=(I-S\,\rchi)^{-1}\,S\,\rchi=(I-S\,\rchi)^{-1}-I$
(the resolvent operator of $S\,\rchi$), whose kernel we denote by
$R(x,y)$, and
$Q_{\epsilon}:=(I-S\,\rchi)^{-1}\,\epsilon\,\varphi$. Then
\eqref{operator} factors into

\begin{equation*}A =(I - S\,\rchi)\,B.\end{equation*}
where $B$ is

\begin{align*}I - & (Q_{\epsilon}\,\otimes\,\rchi\,\psi)
  - \frac{\lambda}{\tilde{\lambda}}\,(I+R)\,S\,(1-\lambda\,\rchi)\,\epsilon\,\com{\rchi}{D}\\ & -
  \frac{\lambda}{\tilde{\lambda}}\,(Q_{\epsilon}\,\otimes\,\psi)\,(1-\lambda\,\rchi)\,\epsilon\,\com{\rchi}{D},\qquad \lambda\neq 1.\end{align*} Hence
\begin{equation*}D_{1,N}(s,\lambda)
=D_{2,N}(s,\tilde{\lambda})\,\det(B).\end{equation*} In order to
find $\det(B)$ we use the identity
\begin{equation}\epsilon\,\com{\rchi}{D}=\sum_{k=1}^{2m}
(-1)^{k}\,\epsilon_{k}\otimes\delta_{k},\end{equation} where
$\epsilon_{k}$ and $\delta_{k}$ are the functions
$\epsilon(x-a_{k})$ and $\delta(x-a_{k})$ respectively, and the
$a_{k}$ are the endpoints of the (disjoint) intervals considered,
$J=\cup_{k=1}^{m}(a_{2\,k-1},a_{2\,k})$. We also make use of the
fact that
\begin{equation}a\otimes b\cdot c\otimes d= \inprod{b}{c}\cdot a\otimes d\end{equation}
where $\inprod{.}{.}$ is the usual $L^{2}$--inner product.
Therefore

\begin{align*}(Q_{\epsilon}\otimes\psi)\,(1-\lambda\,\rchi)\,\epsilon\,\com{\rchi}{D} &=
\sum_{k=1}^{2m} (-1)^{k}Q_{\epsilon}\otimes\psi\cdot
(1-\lambda\,\rchi)\,\epsilon_{k}\otimes\delta_{k}
\\ &=\sum_{k=1}^{2m}
(-1)^{k}\inprod{\psi}{(1-\lambda\,\rchi)\,\epsilon_{k}}\,Q_{\epsilon}\otimes\,\delta_{k}.
\end{align*}
It follows that

\begin{equation*}\frac{D_{1,N}(s,\lambda)}{D_{2,N}(s,\tilde{\lambda})}\end{equation*}
equals the determinant of
\begin{align*} I  & -
Q_{\epsilon}\otimes\rchi\psi \\ & -
\frac{\lambda}{\tilde{\lambda}}\,\sum_{k=1}^{2m}
(-1)^{k}\left[(S+R\,S)\,(1-\lambda\,\rchi)\,\epsilon_{k}
+\inprod{\psi}{(1-\lambda\,\rchi)\,\epsilon_{k}}\,Q_{\epsilon}\right]\otimes\delta_{k}.
 \end{align*}
We now specialize to the case of one interval $J=(t,\infty)$, so
$m=1$, $a_{1}=t$ and $a_{2}=\infty$. We write
$\epsilon_{t}=\epsilon_{1}$, and $\epsilon_{\infty}=\epsilon_{2}$,
and similarly for $\delta_{k}$. Writing the terms in the summation
and using the facts that
\begin{equation}\epsilon_{\infty}=-\frac{1}{2},\end{equation}
 and
 \begin{equation}(1-\lambda\,\rchi)\,\epsilon_{t}=-\frac{1}{2}\,(1-\lambda\,\rchi)+(1-\lambda\,\rchi)\,\rchi,\end{equation}
 then yields
\begin{align*}
I - Q_{\epsilon}\otimes\rchi\psi  -
\frac{\lambda}{2\tilde{\lambda}}\,
\left[(S+R\,S)\,(1-\lambda\,\rchi)+
\inprod{\psi}{(1-\lambda\,\rchi)}\,Q_{\epsilon}\right]\otimes(\delta_{t}-\delta_{\infty})\\
\qquad   +\frac{\lambda}{\tilde{\lambda}}
\left[(S+R\,S)\,(1-\lambda\,\rchi)\,\rchi+
\inprod{\psi}{(1-\lambda\,\rchi)\,\rchi}\,Q_{\epsilon}\right]\otimes\delta_{t}
\end{align*}
which, to simplify notation, we write as
\begin{align*}
 I - Q_{\epsilon}\otimes\rchi\psi  - \frac{\lambda}{2\tilde{\lambda}}\,
\left[(S+R\,S)\,(1-\lambda\,\rchi)+
a_{1,\lambda}\,Q_{\epsilon}\right]\otimes(\delta_{t}-\delta_{\infty}) \\
\qquad   +\frac{\lambda}{\tilde{\lambda}}
\left[(S+R\,S)\,(1-\lambda\,\rchi)\,\rchi+
\tilde{a}_{1,\lambda}\,Q_{\epsilon}\right]\otimes\delta_{t},
\end{align*}
where
\begin{equation}a_{1,\lambda} = \inprod{\psi}{(1-\lambda\,\rchi)},\qquad \tilde{a}_{1,\lambda}= \inprod{\psi}{(1-\lambda\,\rchi)\,\rchi}.\end{equation}
 Now we can use the formula:
\begin{equation}\det\left(I-\sum_{i=1}^{L}\alpha_{i}\otimes\beta_{i}\right)=\det\left(\delta_{jk}-\inprod{\alpha_{j}}{\beta_{k}}\right)_{1\leq j,k\leq L}\end{equation}
In this case, $L=3$, and
\begin{align}\alpha_{1}=Q_{\epsilon}, &\qquad
 \alpha_{2} =\frac{\lambda}{\tilde{\lambda}}\,
\left[(S+R\,S)\,(1-\lambda\,\rchi)+ a_{1,\lambda}\,Q_{\epsilon}\right]\nonumber, \\
\alpha_{3}= & -\frac{\lambda}{\tilde{\lambda}}
\left[(S+R\,S)\,(1-\lambda\,\rchi)\,\rchi+
\tilde{a}_{1,\lambda}\,Q_{\epsilon}\right],\nonumber \\ &
\beta_{1}=\rchi\psi, \qquad  \beta_{2}=\delta_{t}-\delta_{\infty},
\qquad \beta_{3}=\delta_{t}.\end{align} In order to simplify the
notation, define

\begin{align}
Q(x,\lambda,t)&:=(I-S\,\rchi)^{-1}\,\varphi, &
P(x,\lambda,t)&:=(I-S\,\rchi)^{-1}\,\psi, \nonumber\\
Q_{\epsilon}(x,\lambda,t)&:=(I-S\,\rchi)^{-1}\,\epsilon\,\varphi,
& P_{\epsilon}(x,\lambda,t)&:=(I-S\,\rchi)^{-1}\,\epsilon\,\psi,
\end{align}

\begin{align}
q_{_{N}}&:=Q(t,\lambda,t), & p_{_{N}}&:=P(t,\lambda,t),\nonumber\\
q_{\epsilon}&:=Q_{\epsilon}(t,\lambda,t), & p_{\epsilon}&:=P_{\epsilon}(t,\lambda,t),\nonumber\\
u_{\epsilon}&:=\inprod{Q}{\rchi\,\epsilon\,\varphi}=\inprod{Q_{\epsilon}}{\rchi\,\varphi},
&
v_{\epsilon}&:=\inprod{Q}{\rchi\,\epsilon\,\psi}=\inprod{P_{\epsilon}}{\rchi\,\psi}, \nonumber\\
\tilde{v}_{\epsilon}&:=\inprod{P}{\rchi\,\epsilon\,\varphi}=\inprod{Q_{\epsilon}}{\rchi\,\varphi},
&
w_{\epsilon}&:=\inprod{P}{\rchi\,\epsilon\,\psi}=\inprod{P_{\epsilon}}{\rchi\,\psi},
\label{inproddef1}\end{align}

\begin{align}
\mathcal{P}_{1,\lambda}&:= \int(1-\lambda\,\rchi)\,P\,d\,x, &
\tilde{\mathcal{P}}_{1,\lambda} &:=
\int(1-\lambda\,\rchi)\,\rchi\,P\,d\,x, \nonumber\\
\mathcal{Q}_{1,\lambda} &:= \int(1-\lambda\,\rchi)\,Q\,d\,x, &
\tilde{\mathcal{Q}}_{1,\lambda} &:=
\int(1-\lambda\,\rchi)\,\rchi\,Q\,d\,x,\nonumber\\
\mathcal{R}_{1,\lambda}&:= \int(1-\lambda\,\rchi)\,R(x,t)\,d\,x, &
\tilde{\mathcal{R}}_{1,\lambda} &:=
\int(1-\lambda\,\rchi)\,\rchi\,R(x,t)\,d\,x.
\label{inproddef2}\end{align} Note that all quantities in
\eqref{inproddef1} and \eqref{inproddef2} are functions of $t$
alone. Furthermore, let
\begin{equation} c_{\varphi} =
\epsilon\,\varphi(\infty)=\frac{1}{2}\int_{-\infty}^{\infty}\varphi(x)\,d\,x,
\qquad c_{\psi} =
\epsilon\,\psi(\infty)=\frac{1}{2}\int_{-\infty}^{\infty}\psi(x)\,d\,x.\end{equation}
From \cite{Trac2} we find
\begin{equation}\lim_{N\to\infty}c_{\varphi}=\sqrt{\frac{\tilde{\lambda}}{2}},\qquad \lim_{N\to\infty}c_{\psi}=0,\end{equation}
and at $t=\infty$,

\begin{align*}\mathcal{P}_{1,\lambda}(\infty) = 2\,c_{\psi}, \quad
& \mathcal{Q}_{1,\lambda}(\infty) =2\,c_{\varphi},\quad
\mathcal{R}_{1,\lambda}(\infty) = 0 ,\nonumber\\
\tilde{\mathcal{P}}_{1,\lambda}(\infty) = &
\tilde{\mathcal{Q}}_{1,\lambda}(\infty) =\tilde{
\mathcal{R}}_{1,\lambda}(\infty) = 0.\end{align*} Hence
\begin{align}
\inprod{\alpha_{1}}{\beta_{1}} & = \tilde{v}_{\epsilon}, \quad
\inprod{\alpha_{1}}{\beta_{2}}=q_{\epsilon}-c_{\varphi}, \quad
\inprod{\alpha_{1}}{\beta_{3}}=q_{\epsilon},
\\
\inprod{\alpha_{2}}{\beta_{1}} &
=\frac{\lambda}{2\,\tilde{\lambda}}
\,\left[\mathcal{P}_{1,\lambda}-a_{1,\lambda}\,(1-\tilde{v}_{\epsilon})\right], \\
\inprod{\alpha_{2}}{\beta_{2}} & =
\frac{\lambda}{2\,\tilde{\lambda}}\,\left[\mathcal{R}_{1,\lambda}
+ a_{1,\lambda}\,(q_{\epsilon}-c_{\varphi})\right]\label{example}, \\
\inprod{\alpha_{2}}{\beta_{3}} & =
\frac{\lambda}{2\,\tilde{\lambda}}\,\left[\mathcal{R}_{1,\lambda}
+ a_{1,\lambda}\,q_{\epsilon}\right], \\
\inprod{\alpha_{3}}{\beta_{1}} & =
-\frac{\lambda}{\tilde{\lambda}}\,\left[\tilde{\mathcal{P}}_{1,\lambda}
- \tilde{a}_{1,\lambda}\,(1-\tilde{v}_{\epsilon})\right], \\
\inprod{\alpha_{3}}{\beta_{2}} & =
-\frac{\lambda}{\tilde{\lambda}}\,\left[\tilde{\mathcal{R}}_{1,\lambda}+\tilde{a}_{1,\lambda}\,(q_{\epsilon}
- c_{\varphi})\right], \\
 \inprod{\alpha_{3}}{\beta_{3}} & =
-\frac{\lambda}{\tilde{\lambda}}\,
\left[\tilde{\mathcal{R}}_{1,\lambda}+\tilde{a}_{1,\lambda}\,q_{\epsilon}\right].
\end{align}
As an illustration, let us do the computation that led to
\eqref{example} in detail. As in \cite{Trac2}, we use the facts
that $S^{t}=S$, and $(S+S\,R^{t})\,\rchi=R$ which can be easily
seen by writing $R=\sum_{k=1}^{\infty}(S\,\rchi)^{k}$. Furthermore
we write $R(x,a_{k})$ to mean
\begin{equation*} \lim_{\substack{y\to a_{k}\\y\in J}}R(x,y).\end{equation*}
In general, since all evaluations are done by taking the limits
from within $J$, we can use the identity
$\rchi\,\delta_{k}=\delta_{k}$ inside the inner products. Thus
\begin{align*}
\inprod{\alpha_{2}}{\beta_{2}} & =
\frac{\lambda}{\tilde{\lambda}}\,\left[
\inprod{(S+R\,S)\,(1-\lambda\,\rchi)}{\delta_{t}-\delta_{\infty}}
+
a_{1,\lambda}\,\inprod{Q_{\epsilon}}{\delta_{t}-\delta_{\infty}}\right]\\
&=\frac{\lambda}{\tilde{\lambda}}\,\left[\inprod{(1-\lambda\,\rchi)}{(S+R^{t}\,S)\,\left(\delta_{t}-\delta_{\infty}\right)}
+ a_{1,\lambda}\left(Q_{\epsilon}(t)-Q_{\epsilon}(\infty)\right)\right]\\
&
=\frac{\lambda}{\tilde{\lambda}}\,\left[\inprod{(1-\lambda\,\rchi)}{(S+R^{t}\,S)\,\rchi\,\left(\delta_{t}-\delta_{\infty}\right)}
+ a_{1,\lambda}\left(q_{\epsilon}-c_{\varphi}\right)\right]\\
&
=\frac{\lambda}{\tilde{\lambda}}\,\left[\inprod{(1-\lambda\,\rchi)}{R(x,t)-R(x,\infty)}
+ a_{1,\lambda}\left(q_{\epsilon}-c_{\varphi}\right)\right]\\
&
=\frac{\lambda}{\tilde{\lambda}}\,\left[\mathcal{R}_{1,\lambda}(t)-\mathcal{R}_{1,\lambda}(\infty)
+ a_{1,\lambda}\left(q_{\epsilon}-c_{\varphi}\right)\right]\\
&
=\frac{\lambda}{\tilde{\lambda}}\,\left[\mathcal{R}_{1,\lambda}(t)
+ a_{1,\lambda}\left(q_{\epsilon}-c_{\varphi}\right)\right].
\end{align*}
We want the limit of the determinant
\begin{equation}\det\left(\delta_{jk}-\inprod{\alpha_{j}}{\beta_{k}}\right)_{1\leq
j,k\leq L},\end{equation} as $N\to \infty$. In order to get our
hands on the limits of the individual terms involved in the
determinant, we will find differential equations for them first as
in \cite{Trac2}.

Row operation on the matrix show that $a_{1,\lambda}$ and
$\tilde{a}_{1,\lambda}$ fall out of the determinant; to see this
add $\lambda\,a_{1,\lambda}/(2\,\tilde{\lambda})$ times row 1 to
row 2 and $\lambda\,\tilde{a}_{1,\lambda}/\tilde{\lambda}$ times
row 1 to row 3. So we will not need to find differential equations
for them. Our determinant is
{\large\begin{equation}\det\left(\begin{array}{ccc} 1 -
\tilde{v}_{\epsilon} & -(q_{\epsilon}-c_{\varphi}) & -q_{\epsilon}
\\[6pt]
 -\frac{\lambda\,\mathcal{P}_{1,\lambda}}{2\,\tilde{\lambda}} & 1 - \frac{\lambda\,\mathcal{R}_{1,\lambda}}{2\,\tilde{\lambda}}
 & -\frac{\lambda\,\mathcal{R}_{1,\lambda}}{2\,\tilde{\lambda}}
\\[6pt] \frac{\lambda\,\tilde{\mathcal{P}}_{1,\lambda}}{\tilde{\lambda}} & \frac{\lambda\,\tilde{\mathcal{R}}_{1,\lambda}}{\tilde{\lambda}}
& 1 +
\frac{\lambda\,\tilde{\mathcal{R}}_{1,\lambda}}{\tilde{\lambda}}
\end{array}\right).\end{equation}}
Proceeding as in \cite{Trac2} we find the following differential
equations
\begin{align} \frac{d}{d\,t}\,u_{\epsilon} & = q_{_{N}}\,q_{\epsilon}, & \frac{d}{d\,t}\,q_{\epsilon} & = q_{_{N}} -q_{_{N}}\, \tilde{v}_{\epsilon} - p_{_{N}}\,u_{\epsilon},\\
 \frac{d}{d\,t}\mathcal{Q}_{1,\lambda} & =
q_{_{N}}\left(\lambda-\mathcal{R}_{1,\lambda}\right), &
\frac{d}{d\,t}\mathcal{P}_{1,\lambda} & =
p_{_{N}}\left(\lambda-\mathcal{R}_{1,\lambda}\right),\label{ex2}\\
 \frac{d}{d\,t}\mathcal{R}_{1,\lambda} & =
-p_{_{N}}\,\mathcal{Q}_{1,\lambda}-q_{_{N}}\,\mathcal{P}_{1,\lambda},
  & \frac{d}{d\,t}\tilde{\mathcal{R}}_{1,\lambda} & =
-p_{_{N}}\,\tilde{\mathcal{Q}}_{1,\lambda}-q_{_{N}}\,\tilde{\mathcal{P}}_{1,\lambda},
\\  \frac{d}{d\,t}\tilde{\mathcal{Q}}_{1,\lambda} & =
q_{_{N}}\left(\lambda-1-\tilde{\mathcal{R}}_{1,\lambda}\right),
 & \frac{d}{d\,t}\tilde{\mathcal{P}}_{1,\lambda}  & =
p_{_{N}}\left(\lambda-1-\tilde{\mathcal{R}}_{1,\lambda}\right).\end{align}
Let us derive the first equation in \eqref{ex2} for example. From
\cite{Trac3} (equation $2.17$), we have
\begin{equation*}
\frac{\partial Q}{\partial t}=-R(x,t)\,q_{_{N}}.
\end{equation*}
Therefore
\begin{align*}
\frac{\partial \mathcal{Q}_{1,\lambda}}{\partial t} & =
\frac{d}{d\,t}\left[\int_{-\infty}^{t}Q(x,t)\,d\,x-(1-\lambda)\,\int_{\infty}^{t}Q(x,t)\,d\,x\right]\\
& = q_{_{N}} + \int_{-\infty}^{t}\frac{\partial Q}{\partial
t}\,d\,x -
(1-\lambda)\left[q_{_{N}}+\int_{\infty}^{t}\frac{\partial
Q}{\partial t}\,d\,x\right]
\\
& = q_{_{N}} - q_{_{N}}\int_{-\infty}^{t}R(x,t)\,d\,x -
(1-\lambda)\,q_{_{N}}+(1-\lambda)\,q_{_{N}}\,\int_{\infty}^{t}R(x,t)\,d\,x \\
& = \lambda\,q_{_{N}}-q_{_{N}}\,\int_{-\infty}^{\infty}(1-\lambda)\,R(x,t)\,d\,x\\
& = \lambda\,q_{_{N}}-q_{_{N}}\,\mathcal{R}_{1,\lambda} =
q_{_{N}}\left(\lambda-\mathcal{R}_{1,\lambda}\right).
\end{align*}
 Now we change variable from $t$ to $s$ where $t=\tau(s)=
2\,\sigma\,\sqrt{N}+\frac{\sigma\,s}{N^{1/6}}$. Then we take the
limit $N\to \infty$, denoting the limits of $ q_{\epsilon},
\mathcal{P}_{1,\lambda},
\mathcal{Q}_{1,\lambda},\mathcal{R}_{1,\lambda},
\tilde{\mathcal{P}}_{1,\lambda} ,\tilde{\mathcal{Q}}_{1,\lambda} ,
\tilde{\mathcal{R}}_{1,\lambda}$ and the common limit of
$u_{\epsilon}$ and $\tilde{v}_{\epsilon}$ respectively by
$\overline{q}, \overline{\mathcal{P}}_{1,\lambda},
\overline{\mathcal{Q}}_{1,\lambda},
\overline{\mathcal{R}}_{1,\lambda},
\overline{\overline{\mathcal{P}}}_{1,\lambda} ,
\overline{\overline{\mathcal{Q}}}_{1,\lambda} ,
\overline{\overline{\mathcal{R}}}_{1,\lambda}$ and $\overline{u}$.
We eliminate $ \overline{\mathcal{Q}}_{1,\lambda}$ and
$\overline{\overline{\mathcal{Q}}}_{1,\lambda}$ by using the facts
that
$\overline{\mathcal{Q}}_{1,\lambda}=\overline{\mathcal{P}}_{1,\lambda}+\lambda\,\sqrt{2}$
and $\overline{\overline{\mathcal{Q}}}_{1,\lambda}=
\overline{\overline{\mathcal{P}}}_{1,\lambda}$. These limits hold
uniformly for bounded $s$ so we can interchange $\lim$ and
$\frac{d}{d\,s}$. Also
$\lim_{N\to\infty}N^{-1/6}q_{_{N}}=\lim_{N\to\infty}N^{-1/6}p_{_{N}}=q
$ , where $q$ is as in \eqref{D2}. We obtain the systems
\begin{equation}\frac{d}{d\,s}\,\overline{u} = -\frac{1}{\sqrt{2}}\,q\,\overline{q}
,
 \qquad
\frac{d}{d\,s}\,\overline{q} =
\frac{1}{\sqrt{2}}\,q\,\left(1-2\,\overline{u}\right),
\end{equation}
\begin{equation}\frac{d}{d\,s}\overline{\mathcal{P}}_{1,\lambda} =
-\frac{1}{\sqrt{2}}\,q\,\left(\overline{\mathcal{R}}_{1,\lambda}-\lambda\right),
\qquad \frac{d}{d\,s}\overline{\mathcal{R}}_{1,\lambda} =
-\frac{1}{\sqrt{2}}\,q\,\left(2\,\overline{\mathcal{P}}_{1,\lambda}+\sqrt{2\tilde{\lambda}}\right),
\end{equation}
\begin{equation}\frac{d}{d\,s}\overline{\overline{\mathcal{P}}}_{1,\lambda} =
\frac{1}{\sqrt{2}}\,q\,\left(1-\lambda-\overline{\overline{\mathcal{R}}}_{1,\lambda}\right),
\qquad
 \frac{d}{d\,s}\overline{\overline{\mathcal{R}}}_{1,\lambda} =
 -q\,\sqrt{2}\,\overline{\overline{\mathcal{P}}}_{1,\lambda}.
\end{equation}
 The change of variables $q\to\mu=\int_{s}^{\infty} q(x)\,d\,x$
transforms these systems into

\begin{equation}\frac{d}{d\,\mu}\overline{u} =
\frac{1}{\sqrt{2}}\,\overline{q}, \qquad
\frac{d}{d\,\mu}\overline{q} =
-\frac{1}{\sqrt{2}}\,\left(1-2\,\overline{u}\right),
\end{equation}
\begin{equation}\frac{d}{d\,\mu}\overline{\mathcal{P}}_{1,\lambda} =
\frac{1}{\sqrt{2}}\,\left(\overline{\mathcal{R}}_{1,\lambda}-\lambda\right),
\qquad \frac{d}{d\,\mu}\overline{\mathcal{R}}_{1,\lambda} =
\frac{1}{\sqrt{2}}\,\left(2\,\overline{\mathcal{P}}_{1,\lambda}+\sqrt{2\tilde{\lambda}}\right),
\end{equation}
\begin{equation}\frac{d}{d\,\mu}\overline{\overline{\mathcal{P}}}_{1,\lambda} =
-\frac{1}{\sqrt{2}}\,\left(1-\lambda-\overline{\overline{\mathcal{R}}}_{1,\lambda}\right),
\qquad
 \frac{d}{d\,\mu}\overline{\overline{\mathcal{R}}}_{1,\lambda} =
 \sqrt{2}\,\overline{\overline{\mathcal{P}}}_{1,\lambda}.
\end{equation} Since $\lim_{s\to \infty}\mu=0$, corresponding to the boundary
values at $t=\infty$ which we found earlier for
$\mathcal{P}_{1,\lambda}, \mathcal{R}_{1,\lambda},
\tilde{\mathcal{P}}_{1,\lambda} ,
\tilde{\mathcal{R}}_{1,\lambda}$, we now have initial values at
$\mu=0$. Therefore
\begin{equation} \overline{\mathcal{P}}_{1,\lambda}(0) = \overline{\mathcal{R}}_{1,\lambda}(0)=
\overline{\overline{\mathcal{P}}}_{1,\lambda}(0)=
\overline{\overline{\mathcal{R}}}_{1,\lambda}(0)=0.\end{equation}
We use this to solve the systems and get
\begin{align}\overline{q} & = \frac{\sqrt{\tilde{\lambda}}-1}{2\,\sqrt{2}}\,e^{\mu} + \frac{\sqrt{\tilde{\lambda}} + 1}{2\,\sqrt{2}}\,e^{-\mu},\\
 \overline{u} & =
\frac{\sqrt{\tilde{\lambda}}-1}{4}\,e^{\mu} -
\frac{\sqrt{\tilde{\lambda}}+1}{4}\,e^{-\mu}+\frac{1}{2},\\
\overline{\mathcal{P}}_{1,\lambda} & =
\frac{\sqrt{\tilde{\lambda}}-\lambda}{2\,\sqrt{2}}\,e^{\mu} +
\frac{\sqrt{\tilde{\lambda}} + \lambda}{2\,\sqrt{2}}\,e^{-\mu} -
\sqrt{\frac{\tilde{\lambda}}{2}},\\
 \overline{\mathcal{R}}_{1,\lambda} & =
\frac{\sqrt{\tilde{\lambda}}-\lambda}{2}\,e^{\mu} -
\frac{\sqrt{\tilde{\lambda}} + \lambda}{2}\,e^{-\mu} +
\lambda, \\
\overline{\overline{\mathcal{P}}}_{1,\lambda} & =
\frac{1-\lambda}{2\,\sqrt{2}}\,(e^{\mu}-e^{-\mu}), \qquad
\overline{\overline{\mathcal{R}}}_{1,\lambda}=\frac{1-\lambda}{2}\,(e^{\mu}+e^{-\mu}-2).
\end{align}
Substituting these expressions into the determinant gives
\eqref{goedet}, namely
\begin{equation} D_{1}(s,\lambda)= D_{2}(s,\tilde{\lambda})\,\frac{\lambda - 1
- \cosh{\mu(s,\tilde{\lambda})} +
\sqrt{\tilde{\lambda}}\,\sinh{\mu(s,\tilde{\lambda})}}{\lambda -
2},\end{equation} where $D_{\beta}=\lim_{N\to\infty}D_{\beta,N}$.

\subsection{The GSE Case}

The GSE case is the easy one. All calculations in \cite{Trac1} and
\cite{Trac2} go through essentially unchanged except for the
trailing factor of $\lambda$. Therefore we will not reproduce them
here.

\subsection{Interlacing property}

The following series of lemmas establish
Corollary~\eqref{interlacingcor}:
\begin{lemma}
Define
\begin{equation}a_{j}=\frac{d^{j}}{d\,\lambda^{j}}\,\sqrt{\frac{\lambda}{2-\lambda}}\,\,\bigg{\vert}_{\lambda=1}.\label{ajdef}\end{equation}
Then $a_{j}$ satisfies the following recursion
\begin{equation} a_{j} = \begin{cases}
          \quad 1 &\text{if} \quad j=0, \\
           \quad (j-1)\,a_{j-1} &\text{for $j\geq 1$, $j$ even,}  \\
           \quad j\,a_{j-1} &\text{for $j\geq 1$, $j$ odd.} \\
         \end{cases}\end{equation}
\label{ajlemma}\end{lemma}
\begin{proof}
Consider the expansion of the generating function
$f(\lambda)=\sqrt{\frac{\lambda}{2-\lambda}}$ around $\lambda=1$

\[f(\lambda)= \sum_{j\geq 0}\frac{a_{j}}{j!}\,(\lambda-1)^{j}=\sum_{j\geq 0} b_{j}\,(\lambda-1)^{j}\]
Since $a_{j}=j!\,b_{j}$, the statement of the lemma reduces to
proving the following recurrence for the $b_{j}$ \begin{equation}
b_{j} =
\begin{cases}
          \quad 1 &\text{if} \quad j=0, \\
           \quad \frac{j-1}{j}\,b_{j-1} &\text{for $j\geq 1$, $j$ even,}  \\
           \quad b_{j-1} &\text{for $j\geq 1$, $j$ odd.} \\
         \end{cases}\label{bjrecurrence} \end{equation}
Let
\[f^{even}(\lambda)=
\frac{1}{2}\left(\sqrt{\frac{\lambda}{2-\lambda}}+\sqrt{\frac{2-\lambda}{\lambda}}\right),
\qquad f^{odd}(\lambda)=
\frac{1}{2}\left(\sqrt{\frac{\lambda}{2-\lambda}}-\sqrt{\frac{2-\lambda}{\lambda}}\right).\]
These are the even and odd parts of $f$ relative to the reflection
$\lambda-1\to-(\lambda-1)$ or $\lambda\to 2-\lambda$. Recurrence
\eqref{bjrecurrence} is equivalent to
\[\frac{d}{d\,\lambda}\,f^{even}(\lambda)=(\lambda-1)\,\frac{d}{d\,\lambda}\,f^{odd}(\lambda)\]
which is easily shown to be true.
\end{proof}
\begin{lemma}\label{flemma}
Define
\begin{equation}f(s,\lambda)=1-\sqrt{\frac{\lambda}{2-\lambda}}\,\,\tanh{\frac{\mu(s,\tilde{\lambda})}{2}},\end{equation}
for $\tilde{\lambda}=2\,\lambda-\lambda^{2}$. Then
\begin{equation}\frac{\partial^{2\,n}}{\partial\,\lambda^{2\,n}}\,f(s,\lambda)\,\,\bigg{\vert}_{\lambda=1}-
\frac{1}{2\,n+1}\,\frac{\partial^{2\,n+1}}{\partial\,\lambda^{2\,n+1}}\,f(s,\lambda)\,\,\bigg{\vert}_{\lambda=1}=
\begin{cases}
          \quad 1 &\text{if $n=0$,}\\
           \quad0 &\text{if $n\geq 1$.}  \\
         \end{cases} \label{flemmaeq}\end{equation}
\end{lemma}
\begin{proof}
The case $n=0$ is readily checked. The main ingredient for the
general case is Fa\'a di Bruno's formula
\begin{equation}
\frac{d^{n}}{d t^{n}}g(h(t))=\sum\frac{n!}{k_{1}!\cdots k_{n}!}
\left(\frac{d^{k}g}{dh^{k}}(h(t))\right)\left(\frac{1}{1!}\frac{dh}{d
t}\right)^{k_{1}}\cdots\left(\frac{1}{n!}\frac{d^{n} h}{d
t^{n}}\right)^{k_{n}}, \label{faa}\end{equation} where
$k=\sum_{i=1}^{n}k_{i}$ and the above sum is over all partitions
of $n$, that is all values of $k_{1},\ldots, k_{n}$ such that
$\sum_{i=1}^{n} i\,k_{i}=n$. We apply Fa\'a di Bruno's formula to
derivatives of the function
$\tanh{\frac{\mu(s,\tilde{\lambda})}{2}}$, which we treat as some
function $g(\tilde{\lambda}(\lambda))$. Notice that for $j\geq 1$,
$\frac{d^{j}\tilde{\lambda}}{d\,\lambda^{j}}\,\,\big{\vert}_{\lambda=1}$
is nonzero only when $j=2$, in which case it equals $-2$. Hence,
in \eqref{faa}, the only term that survives is the one
corresponding to the partition all of whose parts equal $2$. Thus
we have
\begin{align*}\frac{\partial^{2n-k}}{\partial\,\lambda^{2n-k}}\,&\tanh{\frac{\mu(s,\tilde{\lambda})}{2}}\,\,\bigg{\vert}_{\lambda=1}&&\\
 &= \begin{cases}
           0 &\text{if $k=2j+1$, $j\geq 0$}\\
          \frac{(-1)^{n-j}\,(2\,n-k)!}{(n-j)!} \frac{\partial^{n-j}}{\partial\,\tilde{\lambda}^{n-j}}\,\tanh{\frac{\mu(s,\tilde{\lambda})}{2}}\,\,\bigg{\vert}_{\tilde{\lambda}=1} &\text{for $k=2j$, $j\geq 0$}  \\
         \end{cases} \end{align*}

 \begin{align*}\frac{\partial^{2n-k+1}}{\partial\,\lambda^{2n+1-k}}\,&\tanh{\frac{\mu(s,\tilde{\lambda})}{2}}\,\,\bigg{\vert}_{\lambda=1}
 \\
 &= \begin{cases}
           0 &\text{if $k=2j$, $j\geq 0$}\\
         \frac{(-1)^{n-j}\,(2\,n+1-k)!}{(n-j)!} \frac{\partial^{n-j}}{\partial\,\tilde{\lambda}^{n-j}}\,\tanh{\frac{\mu(s,\tilde{\lambda})}{2}}\,\,\bigg{\vert}_{\tilde{\lambda}=1} &\text{for $k=2j+1$, $j\geq 0$}  \\
         \end{cases} \end{align*}
Therefore, recalling the definition of $a_{j}$ in \eqref{ajdef}
and setting $k=2\,j$, we obtain
\begin{eqnarray*}
\frac{\partial^{2\,n}}{\partial\,\lambda^{2\,n}}\,f(s,\lambda)\,\,\bigg{\vert}_{\lambda=1}&=&
\sum_{k=0}^{2\,n}\binom{2\,n}{k}\,\frac{\partial^{k}}{\partial\,\lambda^{k}}\,\sqrt{\frac{\lambda}{2-\lambda}}\,\frac{\partial^{2\,n-k}}{\partial\,\lambda^{2\,n-k}}\,\tanh{\frac{\mu(s,\tilde{\lambda})}{2}}\,\,\bigg{\vert}_{\lambda=1}\\
&=&
\sum_{j=0}^{n}\frac{(2\,n)!\,(-1)^{n-j}}{(2\,j)!\,(n-j)!}\,a_{2\,j}\,\frac{\partial^{n-j}}{\partial\,\tilde{\lambda}^{n-j}}\,\tanh{\frac{\mu(s,\tilde{\lambda})}{2}}\,\,\bigg{\vert}_{\tilde{\lambda}=1}.
\end{eqnarray*}
Similarly, using $k=2\,j+1$ instead yields
\begin{eqnarray*}
\frac{\partial^{2\,n+1}}{\partial\,\lambda^{2\,n+1}}\,f(s,\lambda)\,\,\bigg{\vert}_{\lambda=1}&=&
\sum_{k=0}^{2\,n+1}\binom{2\,n+1}{k}\,\frac{\partial^{k}}{\partial\,\lambda^{k}}\,\sqrt{\frac{\lambda}{2-\lambda}}\,\frac{\partial^{2\,n+1-k}}{\partial\,\lambda^{2\,n+1-k}}\,\tanh{\frac{\mu(s,\tilde{\lambda})}{2}}\,\,\bigg{\vert}_{\lambda=1}\\
&=&
(2\,n+1)\,\sum_{j=0}^{n}\frac{(2\,n)!\,(-1)^{n-j}}{(2\,j)!\,(n-j)!}\,\frac{a_{2\,j+1}}{2\,j+1}\,\frac{\partial^{n-j}}{\partial\,\tilde{\lambda}^{n-j}}\,\tanh{\frac{\mu(s,\tilde{\lambda})}{2}}\,\,\bigg{\vert}_{\tilde{\lambda}=1}\\
&=&
(2\,n+1)\,\frac{\partial^{2\,n}}{\partial\,\lambda^{2\,n}}\,f(s,\lambda)\,\,\bigg{\vert}_{\lambda=1},
\end{eqnarray*}
since $a_{_{2\,j+1}}/(2\,j+1)=a_{_{2j}}$. Rearranging this last
equality leads to \eqref{flemmaeq}.
\end{proof}

\begin{lemma}
Let $D_{1}(s,\lambda)$ and $D_{4}(s,\tilde{\lambda})$ be as in
\eqref{goedet} and \eqref{gsedet}. Then
\begin{equation}D_{1}(s,\lambda)=D_{4}(s,\tilde{\lambda})\,
\left(1-\sqrt{\frac{\lambda}{2-\lambda}}\,\tanh{\frac{\mu(s,\tilde{\lambda})}{2}}\right)^{2}.\end{equation}
\end{lemma}
\begin{proof}
Using the facts that $-1-\cosh{x}=-2\,\cosh^{2}\frac{x}{2}$,
$1=\cosh^{2}{x}-\sinh^{2}{x}$ and
$\sinh{x}=\sinh\frac{x}{2}\,\cosh\frac{x}{2}$ we get

\begin{eqnarray*}
D_{1}(s,\lambda)&=&\frac{-2}{\lambda-2}\,D_{4}(s,\tilde{\lambda})
+ D_{2}(s,\tilde{\lambda})\,\frac{\lambda+
\sqrt{\tilde{\lambda}}\,\sinh{\mu(s,\tilde{\lambda})}}{\lambda-2}\\
&=&\frac{-2}{\lambda-2}\,D_{4}(s,\lambda)+
D_{2}(s,\tilde{\lambda})\,\frac{\lambda\,\cosh^{2}{\frac{\mu(s,\tilde{\lambda})}{2}}+\lambda\,\sinh^{2}{\frac{\mu(s,\tilde{\lambda})}{2}}
+
\sqrt{\tilde{\lambda}}\,\sinh{\mu(s,\tilde{\lambda})}}{\lambda-2}\\
&=& D_{4}(s,\tilde{\lambda}) +
\frac{D_{4}(s,\tilde{\lambda})}{\cosh^{2}\left(\frac{\mu(s,\lambda)}{2}\right)}\,\frac{\lambda\,\sinh^{2}{\frac{\mu(s,\tilde{\lambda})}{2}}
 +
 \sqrt{\tilde{\lambda}}\,\sinh{\mu(s,\tilde{\lambda})}}{\lambda-2}\\
 &=& D_{4}(s,\tilde{\lambda})\,\left(1-\frac{\lambda\,\sinh^{2}{\frac{\mu(s,\tilde{\lambda})}{2}}
 +
 2\,\sqrt{\tilde{\lambda}}\,\sinh\left(\frac{\mu(s,\lambda)}{2}\right)\,\cosh\left(\frac{\mu(s,\lambda)}{2}\right)}{(\lambda-2)\,\cosh^{2}\left(\frac{\mu(s,\lambda)}{2}\right)}\right)\\
 &=& D_{4}(s,\tilde{\lambda})\,\left(1 - 2\,\sqrt{\frac{\lambda}{2-\lambda}}\,\tanh^{2}{\frac{\mu(s,\tilde{\lambda})}{2}}
 + \frac{\lambda}{2-\lambda}\,\tanh^{2}{\frac{\mu(s,\tilde{\lambda})}{2}} \right)\\
&=&
D_{4}(s,\tilde{\lambda})\,\left(1-\sqrt{\frac{\lambda}{2-\lambda}}\,\,\tanh{\frac{\mu(s,\tilde{\lambda})}{2}}\right)^{2}.
\end{eqnarray*}
\end{proof}
For notational convenience, define
$d_{1}(s,\lambda)=D_{1}^{1/2}(s,\lambda),
d_{4}(s,\lambda)=D_{4}^{1/2}(s,\lambda)$. Then

\begin{lemma}For $n\geq 0$,
\begin{equation*}\left[-\frac{1}{(2\,n+1)!}\,\frac{\partial^{2\,n+1}}{\partial\,\lambda^{2\,n+1}}
+\frac{1}{(2\,n)!}\,\frac{\partial^{2\,n}}{\partial\,\lambda^{2\,n}}\right]\,d_{1}(s,\lambda)\,\,\bigg{\vert}_{\lambda=1}=\frac{(-1)^{n}}{n!}\,\frac{\partial^{n}}{\partial\,\lambda^{n}}\,d_{4}(s,\lambda)\,\,\bigg{\vert}_{\lambda=1}.\end{equation*}
\label{dlemma}\end{lemma}
\begin{proof}
Let
\[f(s,\lambda)=1-\sqrt{\frac{\lambda}{2-\lambda}}\,\,\tanh{\frac{\mu(s,\tilde{\lambda})}{2}}\]
by the previous lemma, we need to show that
\begin{equation*}\left[-\frac{1}{(2\,n+1)!}\,\frac{\partial^{2\,n+1}}{\partial\,\lambda^{2\,n+1}}
+\frac{1}{(2\,n)!}\,\frac{\partial^{2\,n}}{\partial\,\lambda^{2\,n}}\right]\,d_{4}(s,\tilde{\lambda})\,f(s,\lambda)\,\,\bigg{\vert}_{\lambda=1}
=\frac{(-1)^{n}}{n!}\,\frac{\partial^{n}}{\partial\,\tilde{\lambda}^{n}}\,d_{4}(s,\tilde{\lambda})\,\,\bigg{\vert}_{\lambda=1}.
\end{equation*} Now formula \eqref{faa} applied to $d_{4}(s,\tilde{\lambda})$
gives
\begin{equation*}\frac{\partial^{k}}{\partial\,\lambda^{k}}\,d_{4}(s,\tilde{\lambda})\,\,\bigg{\vert}_{\lambda=1}=
\begin{cases}
          \quad 0 &\text{if $k=2j+1$, $j\geq 0$,}\\
           \quad\frac{(-1)^{j}\,k!}{j!} \frac{\partial^{j}}{\partial\,\tilde{\lambda}^{j}}\,d_{4}(s,\tilde{\lambda}) &\text{if $k=2j$, $j\geq 0$.}  \\
         \end{cases} \end{equation*}
Therefore
\begin{align*}
-\frac{1}{(2\,n+1)!}\,\frac{\partial^{2\,n+1}}{\partial\,\lambda^{2\,n+1}}&\,d_{4}(s,\tilde{\lambda})\,f(s,\lambda)\,\,\bigg{\vert}_{\lambda=1}\\
&=-\frac{1}{(2\,n+1)!}\,\sum_{k=0}^{2\,n+1}\binom{2\,n+1}{k}\,\frac{\partial^{k}}{\partial\,\lambda^{k}}\,d_{4\,}
\frac{\partial^{2\,n+1-k}}{\partial\,\lambda^{2\,n+1-k}}\,f
\,\,\bigg{\vert}_{\lambda=1}\\
&=
-\sum_{j=0}^{n}\frac{(-1)^{j}}{(2\,n-2\,j+1)!\,j!}\,\frac{\partial^{j}}{\partial\,\tilde{\lambda}^{j}}\,d_{4\,}
\frac{\partial^{2\,n-2\,j+1}}{\partial\,\lambda^{2\,n-2\,j+1}}\,f
\,\,\bigg{\vert}_{\lambda=1}
\end{align*}
Similarly

\begin{align*}
\frac{1}{(2\,n)!}\,\frac{\partial^{2\,n}}{\partial\,\lambda^{2\,n}}\,d_{4}(s,\tilde{\lambda})\,f(s,\lambda)\,\,\bigg{\vert}_{\lambda=1}
&=\frac{1}{(2\,n)!}\,\sum_{k=0}^{2\,n}\binom{2\,n}{k}\,\frac{\partial^{k}}{\partial\,\lambda^{k}}\,d_{4\,}
\frac{\partial^{2\,n-k}}{\partial\,\lambda^{2\,n-k}}\,f
\,\,\bigg{\vert}_{\lambda=1}\\
&=
\sum_{j=0}^{n}\frac{(-1)^{j}}{(2\,n-2\,j)!\,j!}\,\frac{\partial^{j}}{\partial\,\tilde{\lambda}^{j}}\,d_{4\,}
\frac{\partial^{2\,n-2\,j}}{\partial\,\lambda^{2\,n-2\,j}}\,f
\,\,\bigg{\vert}_{\lambda=1}
\end{align*}
Therefore
\begin{align*}
&\left[-\frac{1}{(2\,n+1)!}\,\frac{\partial^{2\,n+1}}{\partial\,\lambda^{2\,n+1}}
+\frac{1}{(2\,n)!}\,\frac{\partial^{2\,n}}{\partial\,\lambda^{2\,n}}\right]\,d_{4}(s,\tilde{\lambda})\,f(s,\lambda)\,\,\bigg{\vert}_{\lambda=1}\\
&\qquad\qquad=
\sum_{j=0}^{n}\frac{(-1)^{j}}{(2\,n-2\,j)!\,j!}\,\frac{\partial^{j}}{\partial\,\tilde{\lambda}^{j}}\,d_{4\,}(s,\tilde{\lambda})
\left[\frac{\partial^{2\,n-2\,j}}{\partial\,\lambda^{2\,n-2\,j}}\,f-
\frac{1}{2\,n-2\,j+1}\,\frac{\partial^{2\,n-2\,j+1}}{\partial\,\lambda^{2\,n-2\,j+1}}\,f\right]
\,\,\bigg{\vert}_{\lambda=1}
\end{align*}
Now Lemma~\ref{flemma} shows that the square bracket inside the
summation is zero unless $j=n$, in which case it is $1$. The
result follows.
\end{proof}
Lemma~\ref{dlemma} establishes the inductive step in the proof of
Corollary~\ref{interlacingcor}.
\section{Numerics}
Let
\begin{equation}q_{n}(x)=\frac{\partial^{n}}{\partial\lambda^{n}}\,q(x,\lambda)\bigg\vert_{\lambda=1},\end{equation}
so that $q_{0}$ equals $q$ from \eqref{pII}.  In order to compute
$F_{\beta}(s,m)$ it is crucial to know $q_{n}$ accurately.
Asymptotic expansions for $q_{n}$ at $-\infty$ are given in
\cite{Trac3}. We outline how to compute $q_{0}$ and $ q_{1}$ as an
illustration. From \cite{Trac3}, we know that, as $t\to+\infty$
\begin{align}
q_{0}(-t/2) & =\frac{1}{2}\sqrt{t}\left(1-\frac{1}{t^{3}}-\frac{73}{2t^{6}}-\frac{10657}{2t^{9}}-\frac{13912277}{8t^{12}}+\textrm{O}\left(\frac{1}{t^{15}}\right)\right), \nonumber \\
q_{1}(-t/2)
 &
 =\frac{\exp{(\frac{1}{3}t^{3/2})}}{2\sqrt{2\pi}\,t^{1/4}}\left(1+\frac{17}{24t^{3/2}}+\frac{1513}{2^{7}3^{2}t^{3}}+\frac{850193}{2^{10}3^{4}t^{9/2}}-\frac{407117521}{2^{15}3^{5}t^{6}}+\textrm{O}\left(\frac{1}{t^{15/2}}\right)\right).
\end{align}
Quantities needed to compute $F_{\beta}(s,m), m=1,2,$ are not only
$q_{0}$ but also integrals involving $q_{0}$, such as
\begin{equation}
I_{0}=\int_{s}^{\infty}(x-s)\,q_{0}^{2}(x)\,d\,x,\quad
J_{0}=\int_{s}^{\infty}q_{0}(x)\,d\,x.
\end{equation}
Instead of computing these integrals afterwards, it is better to
include them as variables in a system together with $q_{0}$, as
suggested in \cite{Pers1}. Therefore all quantities needed are
computed in one step, greatly reducing errors, and taking full
advantage of the powerful numerical tools in MATLAB. Since
\begin{equation}
I_{0}'=-\int_{s}^{\infty}q_{0}^{2}(x)\,d\,x,\quad
I_{0}''=q_{0}^{2},\qquad J_{0}'=-q_{0},
\end{equation}
the system closes, and can be concisely written
\begin{equation}
\frac{d}{ds}\left(\begin{array}{c}q_{0} \\ q_{0}' \\ I_{0} \\ I_{0}' \\
J_{0}
\end{array}\right) =\left(\begin{array}{c}q_{0}' \\ s\,q_{0}+2q_{0}^3 \\ I_{0}' \\ q_{0}^{2} \\ -q_{0}
\end{array}\right).
\label{q0sys}\end{equation} We first use the MATLAB built--in
Runge--Kutta based ODE solver \texttt{ode45} to obtain a first
approximation to  the solution of \eqref{q0sys} between $x=6$, and
$x=-8$, with an initial values obtained using the Airy function on
the right hand side. Note that it is not possible to extend the
range to the left due to the high instability of the solution a
little after $-8$; (This is where the transition region between
the three different regimes in the so--called ``connection
problem'' lies. We circumvent this limitation by patching up our
solution with the asymptotic expansion to the left of $x=-8$.).
The approximation obtained is then used as a trial solution in the
MATLAB boundary value problem solver \texttt{bvp4c}, resulting in
an accurate solution vector between $x=6$ and $x=-10$.

Similarly, if we define
\begin{equation}
I_{1}=\int_{s}^{\infty}(x-s)\,q_{0}(x)\,q_{1}(x)\,d\,x,\quad
J_{1}=\int_{s}^{\infty}q_{0}(x)\,q_{1}(x)\,d\,x,
\end{equation}
then we have the first--order system
\begin{equation}
\frac{d}{ds}\left(\begin{array}{c}q_{1} \\ q_{1}' \\ I_{1} \\ I_{1}' \\
J_{1}
\end{array}\right) =\left(\begin{array}{c}q_{1}' \\ s\,q_{1}+6q_{0}^2\,q_{1} \\ I_{1}' \\ q_{0}\,q_{1} \\ -q_{0}\,q_{1}
\end{array}\right),
\end{equation}
which can be implemented using \texttt{bvp4c} together with a
``seed'' solution obtained in the same way as for $q_{0}$. Work is
in progress to provide publicly downloadable versions of the
MATLAB routines.

 Table~\ref{table} shows a comparison
of percentiles of the $F_{1}$ distribution to corresponding
percentiles of empirical Wishart distributions. Here $\lambda_{i}$
denotes the $i^{th}$ largest eigenvalue in the Wishart Ensemble.
The percentiles in the $\lambda_{i}$ columns were obtained by
finding the ordinates corresponding to the $F_{1}$--percentiles
listed in the first column, and computing the proportion of
eigenvalues lying to the left of that ordinate in the empirical
distributions for the $\lambda_{i}$. The bold entries correspond
to the levels of confidence most commonly used in statistical
applications. The reader should compare this table to a similar
one in \cite{John1}.
\begin{center}
\begin{table}
\begin{tabular}[htb]{c|ccc|ccc|}
\cline{2-7}
& \multicolumn{3}{|c|}{$\mathbf{100\times 100}$} & \multicolumn{3}{c|}{$\mathbf{100\times 400}$} \\
\hline
\multicolumn{1}{|c|}{\textbf{$F_{1}$-Percentile}} & $\lambda_{1}$ & $\lambda_{2}$ & $\lambda_{3}$ & $\lambda_{1}$ & $\lambda_{2}$ & $\lambda_{3}$ \\
\hline
\multicolumn{1}{|c|}{$0.01$} & $0.008$ & $0.005$  & $0.004$  & $0.008$ & $0.006$ & $0.004$ \\
\multicolumn{1}{|c|}{$0.05$} & $0.042$ & $0.033$  & $0.025$  & $0.042$ & $0.037$ & $0.032$ \\
\multicolumn{1}{|c|}{$0.10$} & $0.090$ & $0.073$  & $0.059$  & $0.088$ & $0.081$ & $0.066$ \\
\multicolumn{1}{|c|}{$0.30$} & $0.294$ & $0.268$  & $0.235$  & $0.283$ & $0.267$ & $0.254$ \\
\multicolumn{1}{|c|}{$0.50$} & $0.497$ & $0.477$  & $0.440$  & $0.485$ & $0.471$ & $0.455$ \\
\multicolumn{1}{|c|}{$0.70$} & $0.699$ & $0.690$  & $0.659$  & $0.685$ & $0.679$ & $0.669$ \\
\multicolumn{1}{|c|}{$\mathbf{0.90}$} & $\mathbf{0.902}$  & $\mathbf{0.891}$  & $\mathbf{0.901}$ & $\mathbf{0.898}$ & $\mathbf{0.894}$ & $\mathbf{0.884}$ \\
\multicolumn{1}{|c|}{$\mathbf{0.95}$} & $\mathbf{0.951}$  & $\mathbf{0.948}$  & $\mathbf{0.950}$ & $\mathbf{0.947}$ & $\mathbf{0.950}$ & $\mathbf{0.941}$ \\
\multicolumn{1}{|c|}{$\mathbf{0.99}$} & $\mathbf{0.992}$  & $\mathbf{0.991}$  & $\mathbf{0.991}$ & $\mathbf{0.989}$ & $\mathbf{0.991}$ & $\mathbf{0.989}$ \\
\hline
\end{tabular}
\caption{Percentile comparison of $F_{1}$ vs. empirical
distributions for $100\times 100$ and $100\times 400$ Wishart
matrices with identity covariance.}\label{table}
\end{table}
\end{center}
\section*{Acknowledgements}
The author would like to thank Craig A. Tracy for the discussions
that initiated this work and for the invaluable guidance and
support that helped complete it, as well as Eric Rains for useful
discussions. This work was supported in part by the National
Science Foundation under grant DMS--0304414.

\end{document}